\documentclass[12 pt]{article}
\usepackage{amsfonts}
\usepackage{amsthm}
\usepackage{amsmath}
\usepackage{enumerate}
\usepackage{amscd}
\usepackage[latin2]{inputenc}
\usepackage{t1enc}
\usepackage[mathscr]{eucal}
\usepackage{indentfirst}
\usepackage{graphicx}
\usepackage{graphics}
\usepackage{mathtools}

\DeclarePairedDelimiter\floor{\lfloor}{\rfloor}
\usepackage[section]{placeins}
\usepackage{pict2e}
\usepackage{epic}
\usepackage{color}
\usepackage{caption}
\usepackage{subcaption}
\usepackage[dvipsnames]{xcolor} 
\numberwithin{equation}{section}
\usepackage[margin=2.9cm]{geometry}
\usepackage{epstopdf} 
\usepackage{hyperref} 
\usepackage{adjustbox}

\usepackage{float}
\usepackage{subcaption}

\theoremstyle{plain}
\newtheorem{Th}{Theorem}[section]
\newtheorem{Lemma}[Th]{Lemma}

\theoremstyle{definition}

\newtheorem{Rem}[Th]{Remark}
\newtheorem{?}[Th]{Problem}

\begin{document}
	
	\title{ \textsc{New Approximations for the higher order coefficients in an asymptotic expansion for the Barnes $G$-function}}
	\author{ Aziz Issaka\footnote{Email: aissaka@uncc.edu} \quad  \\  \\ Department of Mathematics and Statistics \\ University of North Carolina at Charlotte \\9201 University City Blvd,\\Charlotte, NC 28223\\
		USA}
	
	\date{\today}

	%
	%
	%

	\maketitle
	
	\begin{abstract} 
		In this paper, we provide new formulas for determining the coefficients appearing in the asymptotic expansion for the Barnes $G$-function as $n$ tends to infinity for certain classes of asymptotic expansion for the Barnes $G$-function. We remark that our formulas can be used to approximate the coefficients appearing in an asymptotic expansion of the ``random matrix factor" from the Keathing-Snaith conjecture and the coefficients appearing in an asymptotic expansion of the ``L\'evy-Khintchine type representation of the reciprocal of the Barnes $G$-function".
	\end{abstract}
	
	
	\section{Introduction and main results}
	\label{sec1}
	Barnes introduced and studied the multiple gamma functions $\Gamma_n(z)$ in a series of papers \cite{B1, B2, B3, B4} published between 1899 and 1904. The multiple gamma functions $\Gamma_n(z)$ is a generalization of the classical Euler gamma function $\Gamma(z)$. A special interest among the multiple gamma functions is the so called double gamma function or Barnes $G$-function and is defined by $G(z)=1/\Gamma_2(z)$. The research on theory of the Barnes $G$-function has wide range of applications in pure mathematics, applied mathematics and theoretical physics. The $G$-function appears in the theory of $p$-adic $L$-functions \cite{ N4, Nogu}, in the study of determinants of the Laplacians on the $n$-dimensional unit sphere \cite{JunC, N4, OPS, QChoi, PSar, VardiI, XW}. The Barnes $G$-function also appears in the study of random matrix theory and in analytic number theory \cite{VAdam,AMYor}. 
	Keating and Snaith \cite{Keating} make the following celebrated conjecture for the moments of the Riemann zeta function in terms of the Barnes $G$-function $G(z)$ and the gamma function $\Gamma(z)$:
		\begin{align}
		\lim\limits_{T\to \infty} \frac{1}{(\log T)^{\lambda^2}}\frac{1}{T}\int_0^T\bigg|\zeta\left(\frac{1}{2} + it \right)\bigg|^{2\lambda}dt=M(\lambda)A(\lambda)
		\end{align}
		where $\lambda$ is any complex number satisfying $\mathcal{R}(\lambda)>-1$, 
		the random matrix factor $M(\lambda)$ 
		\begin{align}
		M(\lambda) =\frac{(G(1 +\lambda ))^2}{G(1+2\gamma)},
		\end{align}
		and the arithmetic factor $A(\lambda)$ given by
		\begin{align}
		A(\lambda) = \prod_{p\in \mathcal{P}} \left[\left(1- \frac{1}{p}\right)^{\lambda^2}\left(\sum_{m=0}^\infty \left(\frac{\Gamma(\lambda  + m)}{m!\Gamma(\lambda)}\right)^2p^{-m}\right)\right],
		\end{align}
		where $\mathcal{P}$ is the set of prime numbers. Nikeghbali and Yor \cite{AMYor} gave a probabilistic interpretation of the Barnes $G$-function and proposed the following ``L\'evy-Khintchine type representation of the reciprocal of the Barnes $G$-function". In particular, Nikeghbali and Yor \cite{AMYor} showed that for any $z\in\mathbb{C}$, such that $\mathcal{R}(z) > -1$, one has
		\begin{align}
		\frac{1}{G(1+z)} = \exp\left\{-\frac{1}{2}(\log (2\pi) -1)z + (1+\gamma)\frac{z^2}{2} + \int_0^\infty \frac{\left(e^{-uz} - 1 +zu -\frac{u^2z^2}{2}\right)}{u(2\sinh(u/2))^2}du\right\}.
		\end{align}
	The $G$-function is an entire function and is known to satisfy the functional equation 
	\begin{align*}
	&G(z+1) = \Gamma(z) G(z), \;\; z\in\mathbb{C},\\
	&G(0)=G(2)=G(3)=1.
	\end{align*}
	The $G$-function has the following Weierstrass canonical product:
	\begin{align*}
	G(z+1) = (2\pi)^{z/2}\exp\left(-z(z+1)/2 -\gamma z^2/2\right)\prod_{k=1}^\infty\left[\left(1 + \frac{z}{k}\right)^k\exp\left((-z+z^2/(2k)\right)\right],
	\end{align*}
	where $\gamma=0.5772166\dots$ is the Euler-Mascheroni constant. 
	
	A deeper study of the Barnes $G$-function in the January of 2000 was proposed by Richard Askey at the panel discussion of the San Diego symposium on asymptotics and applied analysis. Motivated by Askey's proposition and growing interest of the Barnes $G$-function, Ferreira and L\'opez \cite{{FL}} presented the following asymptotic expansion of logarithm of the Barnes $G$-function in terms of the Bernoulli numbers
	\begin{align}
	\label{EQ11}
	\log G(z+1)&\sim \frac{1}{4}z^2 + z\log\Gamma(z+1) - \left(\frac{1}{2}z(z+1) + \frac{1}{12}\right)\log z -\log A\\
	& + \sum_{n\geq1}\frac{B_{2n + 2}}{2n(2n+1)(2n+2)z^{2n}}\nonumber,
	\end{align} 
	as $z\to +\infty$ in the sector $|\arg(z)| < \pi-\delta$ with any fixed $0<\delta \leq \pi$. Here $B_n$ denotes the $n$th Bernoulli number. The constant $A$ first appeared in the articles of Kinkelin \cite{VAdam, Kinkelin} and Glaisher \cite{VAdam, Glaisher} on the asymptotic behavior of the following product:
	\begin{align*}
	1^1 \cdot 2^2 \cdot 3^3 \cdots n^n = \frac{n!^n}{G(n+1)}
	\end{align*}
	as $n\to +\infty$. The constant $A$ is called Galisher--Kinkelin constant and can be defined by
	\begin{align*}
	\log A = \frac{\gamma + \log(2\pi)}{12} -\frac{\zeta'(2)}{2\pi^2}=\frac{1}{12}-\zeta'(-1)=0.24875447\dots,
	\end{align*}
	where $\zeta$ is the Riemann zeta function \cite{N4}.  Adamchik \cite{VAdam, VAdam1} presented some special values of the Barnes $G$-function and initiated some discussions on algorithms for high-precision computation of the Barnes $G$-function. Recently Chen \cite{Chen1} established the following two general asymptotic expansions for the Barnes $G$-function. The first asymptotic expansion is
		\begin{align}
		\label{EQ12}
		G(z+1)\sim A^{-1}z^{-\frac{1}{2}z^2 - \frac{1}{2}z - \frac{1}{12}}e^{\frac{z^2}{4}}\Gamma^z(z+1)\left(1 + \sum_{n\geq1}\frac{b_n(\ell,r)}{z^n} \right)^{z^{\ell}/r},
		\end{align}
		for $z\to +\infty$ and $|\arg(z)|<\pi$, where $r$ is a nonzero real number, $\ell$ is a nonnegative integer and the coefficients $b_n(\ell,r)(n\in \mathbb{N})$ are given by
		\begin{align}
		\label{EQ14}
		b_n(\ell,r) &= \sum_{(1+\ell)k_1+(2+\ell)k_2+\cdots +(n+\ell)k_n=n}\frac{r^{k_1+k_2+\cdots+k_n}}{k_1!k_2!\cdots k_n!}\cdot \left(\frac{B_3}{1\cdot 2\cdot 3}\right)^{k_1}\cdot \left(\frac{B_4}{2\cdot 3\cdot 4}\right)^{k_2} \nonumber\\
		&\cdots \left(\frac{B_{n+2}}{n\cdot (n+1)\cdot (n+2)}\right)^{k_n}.
		\end{align}
		The second asymptotic expansion is
		\begin{align}
		\label{EQ13}
		&G(z+1)\sim A^{-1}z^{-\frac{1}{2}z^2 - \frac{1}{2}z - \frac{1}{12}}e^{\frac{z^2}{4}}\Gamma^z(z+1)\left(1 + \ln\left(1+\sum_{n\geq 1}\frac{a_n(\ell,r)}{z^n}\right) \right)^{z^{\ell}/r},
		\end{align}
		for $z\to +\infty$ and $|\arg(z)|<\pi$, where $r$ is a nonzero real number, $\ell$ is a nonnegative integer and the coefficients $a_n(\ell,r)(n\in \mathbb{N})$  are given by 
		\begin{align}
		\label{EQ15}
		a_n(\ell,r) &= \sum_{k_1+2k_2+\cdots +nk_n=n}\frac{1}{k_1!k_2!\cdots k_n!}b_1^{k_1}(\ell,r)b_2^{k_2}(\ell,r)\cdots b_n^{k_n}(\ell,r).
		\end{align}
	
	Nemes \cite{N3} established the following new asymptotic expansion for the Gamma function:
	\begin{align}
	\label{EqNemes}
	\left(\frac{\Gamma(x)}{\left(\frac{x}{e}\right)^x\sqrt{\frac{2\pi}{x}}}\right)^x \sim \sum_{n\geq 0}\frac{G_n}{x^n},
	\end{align}
	holds as $x\to \infty$. In \cite{N2}, Nemes proved that  \eqref{EqNemes} also holds in the sector $|\arg z|\leq \pi -\delta$, $0<\delta\leq \pi$ as $z\to \infty$ and went further to study the growth of the coefficients $G_n$ as $n$ becomes large. 
	
	Xu and Wang \cite{XW} established a similar asymptotic expansion for the Barnes $G$-function
	\begin{align}
	\label{XuWang}
	\left(\frac{G(z+1)}{A^{-1}z^{-\frac{1}{2}z^2 - \frac{1}{2}z - \frac{1}{12}}e^{\frac{z^2}{4}}\Gamma^z(z+1)}\right)^r\sim \sum_{n\geq0}\frac{w_n(r)}{z^{2n}},
	\end{align}
	for $z\to +\infty$ and $|\arg(z)|<\pi$, where $r$ is a nonzero real number and the sequence $(w_n(r))_{n\geq 1}$ satisfy the recurrence relation
	\begin{align}
	\label{EQ17}
	w_0(r)= 1,\; w_n(r) =\frac{r}{2n}\sum_{k=0}^{n-1}\frac{B_{2n-2k+2}}{(2n-2k+1)(2n-2k+2)}w_k(r),\;\; k\geq 1.
	\end{align}
	A similar asymptotic expansion \eqref{XuWang} with a recurrence formula was initially studied by Xu \cite{AminXu} when $z$ is a natural number. Throughout this paper, if not stated otherwise, empty sums are taken to be zero. In this paper, we denote the Stirling numbers of the first kind by $s(n,\nu)$ and the Stirling numbers of the second kind by $S(\nu,k)$. We define the Stirling numbers of the second kind $S(\nu,k)$ by the generating function \cite{Comtet, N2, PG}):
	\begin{align}
	\label{1.8}
	\frac{x^k}{(1-x)(1-2x)\cdots (1-kx)}=\sum_{\nu\geq k}S(\nu,k)x^\nu,
	\end{align}
	for every nonnegative integer $k$. The Stirling numbers of the first kind $s(n,\nu)$ is defined by the generating function
	\begin{align}\label{EQ119}
	\{(1+ z)^x - 1\}^k: = \sum_{n\geq k}\sum_{\nu=k}^{n}k!S(\nu,k)s(n,\nu)\frac{x^j}{n!} z^n,
	\end{align}
	for $|z|<1$ and where $k\geq 1$ is a natural number. It known that the Stirling numbers of the second kind can be computed explicitly using the following expression \cite{N1}
	\begin{align*}
	S(\nu,k) = \frac{1}{k!}\sum_{j=0}^{k}(-1)^k \binom kj (k-j)^\nu.
	\end{align*}
	Nemes \cite{N1} proved that the Stirling coefficients $\gamma_n$ appearing in the following asymptotic expansion
	\begin{align}
	\label{1.10}
	\nu! \sim \left(\frac{\nu}{e}\right)^\nu \sqrt{2\pi \nu}\sum_{n\geq 0}\frac{\gamma_n}{\nu^n}
	\end{align}
	as $\nu\to +\infty$ can be computed explicitly using following expression
	\begin{align*}
	\gamma_n =\frac{\Gamma(3n +\frac{3}{2})}{\sqrt{\pi}}\sum_{k=0}^{2n}\frac{2^{n+k+1}}{(2n+2k+1)((2n-2k)}\sum_{j=0}^{k}\frac{(-1)^{k-j}S(2n+2k-j,k-j)}{j!(2n+2k-j)!}.
	\end{align*}
	
	Because of the many applications of the Barnes $G$-function, computer algebra researchers are interested in implementing this function in computer algebra systems. The asymptotic expansion \eqref{XuWang} is a special case of the asymptotic expansion \eqref{EQ12}. This means that the complicated formula $b_n(0,r)$ in \eqref{EQ14} can be computed using the seemingly easy recurrence relation $w_n(r)$ appearing in \eqref{EQ17}. It turn out that calculations of higher order terms by the recurrence relation \eqref{EQ17} can be exceedingly labor intensive. However, its often not necessary to know the exact value of the late terms in the theory of asymptotic expansions and accurate numerical approximations of late terms are enough for application purposes \cite{N2}.
	
	The aim of this paper is in two folds. The first is to provide efficient formulas for the computation of the coefficients appearing in \eqref{EQ17}. The second is to describe the asymptotic behavior of the coefficients \eqref{EQ14} and \eqref{EQ15} for large $n$ when $r$ is fixed. In our first theorem, we provide the following formula for computing the coefficients appearing in \eqref{XuWang}. Henceforth we denote these coefficients by $A_n(r)$.

	\begin{Th}
		\label{thm111} Let $r\neq 0$ be a real number. Then we have
	    \begin{align}
	    \label{EQQ112}
	    G(z+1)\sim A^{-1}z^{-\frac{1}{2}z^2 - \frac{1}{2}z - \frac{1}{12}}e^{\frac{z^2}{4}}\Gamma^z(z+1)\left(1 + \sum_{n\geq1}\frac{A_n(r)}{z^n} \right)^{1/r},
	    \end{align}
		for $z\to\infty$ and $|\arg(z)|<\pi$, where the sequence $(A_n(r))_{n\geq 0}$ can be determined by 
		\begin{equation}
		\label{EQ111}
		A_n(r)= \sum_{k=0}^{m-1}\frac{rB_{2n-2k+2}A_k(r)}{(2n-2k)(2n-2k+1)(2n-2k+2)} + \mathcal{O}\left(\frac{rB_{2n-2m+2}}{(2n-2m)(2n-2m+1)(2n-2m+2)} \right),
		\end{equation}
		for every integer $m>0$ as $n$ tends to infinity.
	\end{Th}
     In our second theorem, we give an asymptotic formula for the coefficients $A_n(r)$ when $n$ is large and $r$ is fixed.
	\begin{Th}
		\label{thm112a} Let $r\neq 0$ be a real number. Then we have
		\begin{align}
		\label{EQQ112b}
		G(z+1)\sim A^{-1}z^{-\frac{1}{2}z^2 - \frac{1}{2}z - \frac{1}{12}}e^{\frac{z^2}{4}}\Gamma^z(z+1)\left(1 + \sum_{n\geq1}\frac{A_n(r)}{z^n} \right)^{1/r},
		\end{align}
		for $z\to\infty$ and $|\arg(z)|<\pi$, where the coefficients $A_n(r)$ have the following asymptotic expansion
		\begin{align}
		\label{thm112}
		A_n(r) \sim (-1)^{n}r(2n-1)\left(\frac{n-1}{\pi e}\right)^{2n-2}\sqrt{\frac{n-1}{\pi^7}}\sum_{k\geq 0}\frac{I_{k,n}(r)}{(n-1)^k},
		\end{align}
		as $n\to +\infty$, where $I_{0,n}= 1/4$ and
		\begin{align*}
		I_{k,n}(r)= \frac{\gamma_k}{2^{k+2}} + \frac{1}{(2n-1)}\sum_{j=1}^{k}\sum_{\nu=1}^{\floor{\frac{j+1}{2}}}(-1)^\nu2^{2\nu-k-2}\pi^{2\nu}A_\nu(r)\gamma_{k-j}S(j-1, 2\nu-2),
		\end{align*}
		for $k\geq 1$.
	\end{Th}
	\begin{Rem} It is worth mentioning that the recurrence relation \eqref{EQ17} proposed by Xu and Wang \cite{XW} quickly runs into a great difficulty when computing $w_n(r)$ for even small values of $n$ and fixed $r$. Using Wolfram Mathematica 12 software with our formulas \eqref{EQ111} and \eqref{thm112}, we are able to provide approximate values for these coefficients. Our formulas can be applied to the results of Xu \cite{AminXu}. Because of importance of the Keating-Snaith conjecture, it seems interesting to provide an asymptotic expansion for the ``random matrix factor" $M(\lambda)$. We present the following asymptotic expansion for the ``random matrix factor" $M(\lambda)$:
			\begin{align}
			M(\lambda)\sim A^{-1}2^{2\lambda^2+\lambda+\frac{1}{12}}z^{\lambda^2 - \frac{1}{12}}\frac{\Gamma^{2\lambda}(\lambda+1)}{\Gamma^{2\lambda}(2\lambda+1)}e^{-\frac{\lambda^2}{2}}\sum_{n\geq0}\frac{J_n}{\lambda^{2n}},
			\end{align}
			for $\lambda\to\infty$ and $|\arg(\lambda)|<\pi$, where the coefficients $J_n$ are given by
			\begin{align*}
			J_n = \sum_{\ell=0}^n\frac{A_{\ell}(2)A_{n-\ell}(-1)}{4^{n-\ell}}.
			\end{align*}
			 Our formulas \eqref{EQ111} and \eqref{thm112} can be used to approximate these $J_n$ coefficients as well as the following coefficients appearing in the ``L\'evy-Khintchine type representation of the reciprocal of the Barnes $G$-function":
			\begin{align}
			\frac{1}{G(z+1)}\sim Az^{\frac{1}{2}z^2 + \frac{1}{2}z+\frac{1}{12}}e^{-\frac{z^2}{4}}\Gamma^{-z}(z+1)\sum_{n\geq 0}\frac{A_n(-1)}{z^{2n}},
			\end{align}
			for $z\to\infty$ and $|\arg(z)|<\pi$.
	\end{Rem}
	\subsection{The proof of Theorem \ref{thm111} and Theorem \ref{thm112a}}
	To set the stage, we need several lemmas. The first two auxilliary theorems are about the asymptotics of the coefficients of certain formal power series.
	\begin{Th}[E. A. Bender, \cite{EA}]\label{thm21} Suppose that 
		\begin{equation*}
		\alpha(x) = \sum_{n\geq 1}\alpha_{n}x^n, \;    \;    F(x,y)= \sum_{i,j\geq 0}f_{ij}x^iy^j,  
		\end{equation*}
		\begin{equation*}
		\beta(x)= F(x,\alpha(x)):= \sum_{n\geq 0}\beta_{n}x^n, \;    \; D(x) =\sum_{n\geq 0}d_nx^n = \frac{\partial F(x,y)}{\partial y}\Bigr|_{\substack{y=\alpha(x)}}.
		\end{equation*}
		Assume that $F(x,y)$ is analytic in $x$ and $y$ in a neighborhood of $(0,0)$, $\alpha_n \neq 0$ and 
		\begin{enumerate}[\normalfont(i)]
			\item \label{cond 1}$\begin{aligned}[t] \alpha_{n-1} =o(\alpha_n)  \;  \text{as} \; n\to +\infty   \end{aligned}$,
			\item \label{cond 2}$\begin{aligned} \sum_{k=m}^{n-m}|\alpha_k\alpha_{n-k}| =\mathcal{O}(\alpha_{n-m}) \;\; \text{for some}\;\; m>0 \;\;  \text{as} \; n\to +\infty  \end{aligned}$, 
			\item \label{cond 3}$|f_{ij}\alpha_{n-i-j+1}| \leq K(\delta)\delta^{i+j}|\alpha_{n-m}| $ when $n\geq i+j>m$ and  $\delta>0$ for some $K(\delta)$.
		\end{enumerate}
		Then 
		\begin{equation*}
		\beta_n =\sum_{k=0}^{m-1} d_k\alpha_{n-k} + \mathcal{O}(\alpha_{n-m}) \;\; \text{as} \;\; n \to +\infty.
		\end{equation*}
	\end{Th} 
	We remark that the assumptions on $F$ are automatically satisfied if $F$ is analytic at the origin (see Odlyzko \cite[p. 116]{Odlyzko}). 
	\begin{Th}[E.A. Bender and L.B. Richmond, \cite{EAL}]\label{thm22} Let $\alpha(x)$ be a formal power series with coefficients $\alpha_k$ where $\alpha_0\neq 0$. Let $p_n$ be the coefficient of $x^n$ in $(1+ \alpha(x))^{\theta n + \kappa}$ where $\theta\neq 0$ and $\kappa$ are fixed complex numbers. If 
		\begin{enumerate}[\normalfont(i)]
			\item $n\alpha_{n-1} =o(\alpha_n)$ as $n\to +\infty$,
		\end{enumerate}
		then
		$p_n\sim(\theta n + \kappa)\alpha_n$ as $n\to +\infty$.
	\end{Th}
	The following two lemmas are easy consequences of the well-known formula \cite{AI}
	\begin{align*}
	B_{2n+2} = (-1)^{k}\frac{2(2n+2)!}{(2\pi)^{2n+2}} \left(1+ \frac{1}{2^{2n+2}} + \frac{1}{3^{2n+2}} + \cdots\right).
	\end{align*}
	\begin{Lemma}\label{lem23}
		For any $n\geq 0$, we have that
		\begin{align*}
		\frac{2(2n+2)!}{(2\pi)^{2n}} <|B_{2n+2}| < \frac{4(2n+2)!}{(2\pi)^{2n+2}}. 
		\end{align*}
	\end{Lemma}
	\begin{Lemma}\label{lem24}We have
		\begin{align*}
		(1 - 2^{-1 - 2n})B_{2n+2} \sim B_{2n+2} \sim (-1)^{n}\frac{2(2n+2)!}{(2\pi)^{2n+2}},
		\end{align*}
		as $n \to +\infty$.
	\end{Lemma}
	
	\subsubsection{Proof of Theorem \ref{thm111}} After a simple algebraic manipulations of \eqref{EQ11}, we obtain
	\begin{align*}
	\left(\frac{G(z+1)}{A^{-1}z^{-\frac{1}{2}z^2 - \frac{1}{2}z - \frac{1}{12}}e^{\frac{z^2}{4}}\Gamma^z(z+1)} \right)^r\sim \exp\left(\sum_{n\geq1}^{\infty}\frac{rB_{2n + 2}}{2n(2n+1)(2n+2)z^{2n}}\right)
	\end{align*}
	which holds in the sector $|\arg(z)|<\pi-\delta$, $0<\delta \leq \pi$ as  $z\to +\infty$. Using \eqref{EQ11} we have the following formal expansion
	\begin{equation*}
	\exp\left(\sum_{n\geq1}\frac{rB_{2n+2}}{2n(2n+1)(2n+2)z^{2n}}\right)\sim \sum_{n\geq 0}\frac{A_n(r)}{z^{2n}}.
	\end{equation*}
	From this we have the following formal generating function 
	\begin{equation*}
	\exp\left(\sum_{n\geq1}\frac{rB_{2n+2}}{2n(2n+1)(2n+2)}x^n\right)= \sum_{n\geq 0}A_n(r)x^n.
	\end{equation*}
	We apply Theorem \ref{thm21} to the formal power series
	\begin{equation*}
	\alpha(x) = \sum_{n\geq1}\frac{rB_{2n+2}}{2n(2n+1)(2n+2)}x^n:= \sum_{n\geq1}\alpha_n(x)x^n,
	\end{equation*}
	and
	\begin{equation*}
	F(x,y)=e^y =\sum_{n\geq 0}\frac{y^n}{n!}.
	\end{equation*}
	This implies that we have
	\begin{equation*}
	B(x)=D(x)=\sum_{n\geq 0}A_n(r)x^n.
	\end{equation*}
	Applying Lemma \ref{lem24} to the sequence
	\begin{equation*}
	\alpha_n(r) = \frac{rB_{2n+2}}{2n(2n+1)(2n+2)},
	\end{equation*}
	we obtain
	\begin{equation*}
	\frac{2|r|(2n-1)!}{(2\pi)^{2n+2}}<|\alpha_n(r)|<\frac{4|r|(2n-1)!}{(2\pi)^{2n+2}},
	\end{equation*}
	for every $n\geq 1$. Since $\alpha_n(r)\neq0$ and
	\begin{align*}
	0\leq \lim_{n\to +\infty}\Big|\frac{\alpha_{n-1}(r)}{\alpha_n(r)}\Big|<\lim_{n\to+\infty}\frac{4\pi^2}{(2n-1)(2n-2)}=0,
	\end{align*}
	we have
	\begin{align*}
	\alpha_{n-1}(r) = o(\alpha_n(r)) \;\;\text{as}\;\; n\to +\infty. 
	\end{align*}
	Hence condition $(1)$ of Theorem \ref{thm21} holds. Also condition $(2)$ of Theorem \ref{thm21} holds for every integer $m>0$ since for $n>2m$ 
	
	\begin{align*}
	\sum_{k=m}^{n-m}|\alpha_k(r)\alpha_{n-k}(r)| &<\sum_{k=m}^{n-m}\frac{4|r|(2k-1)!}{(2\pi)^{2k+2}}\frac{4|r|(2(n-k)-1)!}{(2\pi)^{2(n-k)+2}}\\
	&= \frac{2r^2(2n-2m-1)!}{(2\pi)^{2n-2m+2}}\sum_{k=m}^{n-m}\frac{4(2k-1)!}{(2\pi)^{2k+2}}\frac{4(2(n-k)-1)!}{(2\pi)^{2(n-k)+2}}\frac{(2\pi)^{2n-2m+2}}{2(2n-2m-1)!}\\
	&< |\alpha_{n-m}(r)|\frac{8|r|}{(2\pi)^{2m+2}(2n-2m-1)!}\sum_{k=m}^{n-m}(2k-1)!(2n-2k-1)!\\
	&<|\alpha_{n-m}(r)|\frac{8}{(2\pi)^{2r+2}(2n-2m-1)!}(2(2n-2m-1)!(2m-1)! \\
	&+ (n-2m-1)(2m+1)!(2n-2m-3)!)\\
	&<|\alpha_{n-m}(r)|\frac{8|r|}{(2\pi)^{2m+2}(2n-2m-1)!}(2(2m-1)! + (2m+1)!).
	\end{align*}
	Because $F(x,y) = e^y$ is analytic in $x$ and $y$, it follows that
	\begin{equation*}
	A_n(r) = \sum_{k=0}^{m-1}A_k(r)\frac{rB_{2n-2k+2}}{(2n-2k)(2n-2k+1)(2n-2k+2)} + \mathcal{O}\left(\frac{rB_{2n-2m+2}}{(2n-2m)(2n-2m+1)(2n-2m+2)} \right)
	\end{equation*} 
	for every integer positive integer $m$ as $n\to +\infty$.
	\subsubsection{Proof of Theorem \ref{thm112a}} Substituting
	\begin{equation*}
	\frac{rB_{2n+2}}{2n(2n+1)(2n+2)} = (-1)^n\frac{2r(2n-1)!}{(2\pi)^{2n+2}}\zeta(2n+2)
	\end{equation*}
	in \eqref{EQ111}, we obtain
	\begin{align*}
	A_n(r) = (-1)^n\frac{2(2n-1)!}{(2\pi)^{2n+2}}\zeta(2n+2)\Bigg( \Bigg.&\sum_{k=0}^{m-1}(-1)^k(2\pi)^{2k}A_k(r)\frac{r(2n+2k-1)!}{(2n-1)!}\frac{\zeta(2n-2k+2)}{\zeta(2n+2)} \\
	&+ \mathcal{O}\left(\frac{r(2n-2m-1)!}{(2n-1)!}\right)\Bigg. \Bigg).
	\end{align*}
	Using the definition of the Riemann zeta function we have $\zeta(s) = 1 + \mathcal{O}(1/2^s)$ for large positive $s$. From this it is easy to see that 
	\begin{equation*}
	\frac{\zeta(2n-2k+2)}{\zeta(2n+2)} = 1 + \mathcal{O}\left(\frac{1}{4^{n-k+1}}\right),
	\end{equation*}
	as $n\to +\infty$. It follows that
	\begin{equation*}
	A_n(r) =(-1)^n\frac{2(2n-1)!}{(2\pi)^{2n+2}}\zeta(2n+2)\left(\sum_{k=0}^{m-1}(-1)^k(2\pi)^{2k}A_k(r)\frac{r(2n-2k-1)!}{(2n-1)!}+\mathcal{O}\left(\frac{r(2n-2m-1)!}{(2n-1)!}\right) \right)
	\end{equation*}
	for every integer $m>0$ as $n\to +\infty$. We can write this expression as
	\begin{equation*}
		A_n(r) =(-1)^n\frac{2r(2n-1)!}{(2\pi)^{2n+2}}\zeta(2n+2)\left(\sum_{k=0}^{m-1}\frac{(-1)^k(2\pi)^{2k}A_k(r)}{(2n-1)_{2k}}+\mathcal{O}\left(\frac{1}{(2n-1)_{2m}}\right) \right),
	\end{equation*}
	where
	\[
	(x)_i= 
	\begin{cases}
	x(x-1)\cdots (x-i+1),& i\geq 1,\\
	1,              & i=0,
	\end{cases}
	\] is the falling factorial. This leads us to the formal asymptotic series
	\begin{equation}
	\label{EQ21}
	A_n(r) \sim (-1)^n\frac{2r(2n-1)!}{(2\pi)^{2n+2}}\zeta(2n+2)\sum_{k\geq0}\frac{(-1)^k(2\pi)^{2k}A_k(r)}{(2n-1)_{2k}}
	\end{equation}
	which holds as $n\to +\infty$. Using the Stirling formula \eqref{1.10} we obtain
	\begin{equation}
	\label{EQ22}
	(-1)^n\frac{2r(2n-1)!}{(2\pi)^{2n+2}}\zeta(2n+2) \sim (-1)^nr(2n-1)\left(\frac{n-1}{\pi e}\right)^{2n-2}\sqrt{\frac{n-1}{\pi^7}}\sum_{k\geq 0}\frac{\gamma_k}{2^k(n-1)^k}.
	\end{equation}
	Using the generating function of the Stirling numbers of the second kind \eqref{1.8} one can easily show that
	\begin{equation*}
	\frac{1}{(2n-1)_{2k}} =\frac{1}{(2n-1)}\sum_{\nu\geq 2k-1}\frac{S(\nu-1,2k-2)}{2^\nu}\frac{1}{(n-1)^\nu},
	\end{equation*}
	for $k\geq 1$. Hence \eqref{EQ21} simplifies to the following
	\begin{align*}
	\sum_{k\geq0}\frac{(-1)^k(2\pi)^{2k}A_k(r)}{(2n-1)_{2k}}&=1+\sum_{k\geq1}\frac{(-1)^k(2\pi)^{2k}A_k(r)}{(2n-1)}\sum_{\nu\geq 2k-1}\frac{S(\nu-1,2k-2)}{2^\nu}\frac{1}{(n-1)^\nu} \\
	&\sim 1+\sum_{\nu\geq1}\left[\sum_{k=1}^{\floor*{\frac{\nu+1}{2}}}\frac{(-1)^k(2\pi)^{2k}A_k(r)}{(2n-1)}\frac{S(\nu-1,2k-2)}{2^\nu}\right]\frac{1}{(n-1)^\nu}
	\end{align*}
	as $n\to +\infty$. After substituting this expression and \eqref{EQ22} into \eqref{EQ21} and performing the product of the asymptotic series, we have
	\begin{equation*}
	A_n(r) \sim (-1)^{n}r(2n-1)\left(\frac{n-1}{\pi e}\right)^{2n-2}\sqrt{\frac{n-1}{\pi^7}}\sum_{k\geq 0}\frac{I_{k,n}(r)}{(n-1)^k}
	\end{equation*}
	as $n\to +\infty$, where
	\begin{equation*}
	I_{k,n}(r)= \frac{\gamma_k}{2^{k+2}} + \frac{1}{(2n-1)}\sum_{j=1}^{k}\sum_{\nu=1}^{\floor{\frac{j+1}{2}}}(-1)^\nu2^{2\nu-k-2}\pi^{2\nu}A_\nu(r)\gamma_{k-j}S(j-1, 2\nu-2)
	\end{equation*}
	for $k\geq 1$.
	\section{Asymptotic formulas for coefficients appearing in general asymptotic expansion for the Barnes $G$-function}
   In this section, we describe the asymptotic behavior of coefficients \eqref{EQ14} and \eqref{EQ15} when $n$ is large and $r$ is fixed. In our first theorem, we give an asymptotic formula for the coefficients \eqref{EQ13} when $n$ is large and $r$ is fixed.
	\begin{Th}\label{thm114} Let $r\neq 0$ be a given real number and $\ell\geq 0$ be a given integer. Then
		\begin{equation}
		\label{EQ114}
		b_n(\ell,r)\sim (-1)^n\frac{2r(2n-1)!}{(2\pi)^{2n+2}}
		\end{equation}
		as $n\to +\infty$.
	\end{Th}
	In our second theorem, we give an asymptotic formula for the coefficients \eqref{EQ15} when $n$ is large and $r$ is fixed.
	\begin{Th}\label{thm115} Let $r\neq 0$ be a given real number and $\ell\geq 0$ be a given integer. Then
		\begin{equation}
		\label{EQ115}
		a_n(\ell,r)\sim (-1)^n\frac{2r(2n-1)!}{(2\pi)^{2n+2}}
		\end{equation}
		as $n\to +\infty$.
	\end{Th}
	\begin{Rem} Although the general terms are given by a complicated formula, their asymptotic behavior are simple and it also show the divergent character of the series.It is also interesting to note that the asymptotics has the ``factorial divided by power'' form, which is a very common behavior in asymptotic series (see, e.g., \cite{Berry, Boyd, AI}).
	\end{Rem}
	\subsection{The proof of Theorem \ref{thm114} and Theorem \ref{thm115}}
	\subsubsection{Proof of Theorem \ref{thm114}}
	To prove formula \eqref{EQ114}, we apply Theorem \ref{thm21} to the formal power series
	\begin{align*}
	\alpha(x)=\sum_{n\geq1}A_n(1) x^n := \sum_{n\geq 1}\alpha_nx^n. 
	\end{align*}
	Using the general asymptotic expression \eqref{EQ111} and Lemma \ref{lem24}, we obtain
	\begin{align*}
	\alpha_n\sim (-1)^n\frac{2(2n-1)!}{(2\pi)^{2n+2}},
	\end{align*}
	as $n\to +\infty$. This further implies that condition $(1)$ and condition $(2)$ in Theorem \ref{thm21} are satisfied. Using \eqref{EQ119}, we apply Theorem \ref{thm21} to the function
	\begin{align*}
	F(x,y):= (1+y)^{rx^{\ell}}-1,
	\end{align*}
	where $\ell\geq 0$ is an integer and $r\neq0$ is a real number. It follows that
	\begin{align*}
	B(x)= F(x, \alpha(x))=(1+\alpha(x))^{rx^\ell}-1, 
	\end{align*}
	and
	\begin{align*}
	D(x)=\sum_{n\geq 0}d_nx^n = \frac{\partial F(x,y)}{\partial y} |_{y=\alpha(x)}.
	\end{align*}
	We have that 
	\begin{align*}
	c_1 \frac{2(2n-1)!}{(2\pi)^{2n+2}} <|\alpha_n| < c_2\frac{2(2n-1)!}{(2\pi)^{2n+2}},
	\end{align*}
	for some constants $c_1>c_2>0$. It is clear that condition $(1)$ and $(2)$ of Theorem \ref{thm21} hold. Since $\alpha_n\neq 0$ and
	\begin{align*}
	\left|\frac{f_{ij}\alpha_{n-i-j+1}}{\alpha_{n-1}}\right| & <\frac{s(j,i)(2n-2i-2j)!(2\pi)^{2n}}{j!(2\pi)^{2(n-i-j)+4}(2n-4)!}< C\frac{(2\pi)^{2(i+j)}}{(2\pi)^4}, 
	\end{align*}
	for $n > i+ j> 1$ and for some constant $C>0$. This implies that condition $(3)$ of Theorem \ref{thm21} is satisfied with $\delta = 4\pi^2$. Hence we conclude that
	\begin{align*}
	b_n(\ell,r) \sim rA_n(1) \sim  (-1)^n\frac{2r(2n-1)!}{(2\pi)^{2n+2}},
	\end{align*}
	as $n\to +\infty$.
	\subsubsection{Proof of Theorem \ref{thm115}} To prove formula \eqref{EQ115} we consider the following as a formal power series
	\begin{align*}
	\sum_{n\geq0}a_n(\ell,r)x^n=\exp\left(\sum_{n\geq1}b_n(\ell,r)x^n\right).
	\end{align*}
	We apply Theorem \ref{thm21} to the function
	\begin{align*}
	\alpha(x) = \sum_{n\geq 1}b_n(\ell,r)x^n,
	\end{align*}
	and $F(x,y) := e^y$. It follows that
	\begin{align*}
	\beta(x)= F(x,\alpha(x)) =\sum_{n\geq0}a_n(\ell,r)x^n.
	\end{align*}
	Since 
	\begin{align*}
	b_n(\ell,r) \sim (-1)^n\frac{2r(2n-1)!}{(2\pi)^{2n+2}},
	\end{align*}
	as $n\to+\infty$, the conditions of Theorem \ref{thm21} are satisfied. Hence we conclude that
	\begin{align*}
	a_n(\ell,r) \sim f_1 b_n(\ell,r) \sim (-1)^n\frac{2r(2n-1)!}{(2\pi)^{2n+2}},
	\end{align*}
	as $n\to +\infty$.
	\section{Numerical examples and applications} In this section we present numerical examples to support our proposed approximations for the coefficients appearing in the Barnes $G$-function. We perform all the computations using the Wolfram Mathematica $12$ software. From \eqref{EQ111} we have the following approximation
	From \eqref{EQ111} we have the following approximation formula
	\begin{align}\label{EQ32}
	A_n(r) \approx \sum_{k=0}^{m-1}A_k(r)\frac{rB_{2n-2k+2}}{(2n-2k)(2n-2k+1)(2n-2k+2)}.
	\end{align}
	Table $1$ shows the numerical performance of this approximation for $n=25$, $r=1$ and $m=5, 10, 15$.
		\FloatBarrier
	\begin{table}[ht]
		\caption*{\textbf{Table 1:} Approximation for $A_{25}(1)$ with various values of $m$, usuing \eqref{EQ32}.} 
		\label{tab:a}	
		\centering 
		\begin{tabular}{|c|c|} 
			\hline 
			Exact numerical value of $w_{25}(1)$ & $-3.80007230719156835910256254456\cdot10^{21}$\\ [0.5ex] 
			\hline
			Value of $m$ & $5$\\
			Approximation for $A_{25}(1)$ & $-3.80007230718902365313203214270\cdot10^{21}$\\
			Error & \hspace{-0.8in} $6.6964672375170116\cdot10^{-13}$\\
			\hline
			Value of $m$ & $10$ \\
			Approximation for $A_{25}(1)$ & $-3.80007230719156771563759457627\cdot10^{21}$\\
			Error & \hspace{-1in} $1.6932966426732\cdot10^{-16}$\\
			\hline
			Value of $m$ & $15$ \\
			Approximation for $A_{25}(1)$ & $-3.80007230719156866209504749935\cdot10^{21}$\\
			Error & \hspace{-1in} $7.973334727905\cdot10^{-17}$\\
			\hline
		\end{tabular}
	\end{table}
	\FloatBarrier
	
	Another family of approximation comes from \eqref{thm112}
	\begin{align}\label{EQ33}
	A_n(r) \approx (-1)^nr(2n-1)\left(\frac{n-1}{\pi e}\right)^{2n-2}\sqrt{\frac{n-1}{\pi^7}}\sum_{k=0}^m\frac{I_{k,n}(r)}{(n-1)^k}.
	\end{align}
	Table $2$ shows the numerical performance of this approximation for $n=25$, $r=1$ and $m=0, 5, 10, 15$.
			\begin{table}[ht]
		\caption*{\textbf{Table 2:} Approximation for $A_{25}(1)$ with various values of $m$, usuing \eqref{EQ33}.}
		\centering 
		\begin{tabular}{|c|c|} 
			\hline
			Exact numerical value of $w_{25}(1)$ & $-3.80007230719156835910256254456\cdot10^{21}$\\[0.5ex]
			\hline
			Value of $m$ & 0\\
			Approximations for $A_{25}(1)$ & $-3.79339211041501539979669705428\cdot 10^{21}$\\
			Error & \hspace{-0.4in} $0.00175791307020942925009466724$\\
			\hline
			Value of $m$ & $5$\\
			Approximations for $A_{25}(1)$ & $-3.80007230665300408064593768490\cdot 10^{21}$\\
			Error & \hspace{-0.65in} $1.4172474493114293334\cdot10^{-10}$\\
			\hline
			Value of $m$ & $10$\\
			Approximations for $A_{25}(1)$ & $-3.80007230719068729982239592063\cdot10^{21}$\\
			Error & \hspace{-0.8in} $2.3185329355423978\cdot10^{-13}$\\
			\hline
			Value of $m$ & $15$\\
			Approximations for $A_{25}(1)$ & $-3.80007230719155206569994906990\cdot 10^{21}$\\
			Error & $4.28765594345131\cdot10^{-15}$\\
			\hline
		\end{tabular}
		\label{tab:b} 
\end{table}
	\FloatBarrier
Table $3$ shows values of our approximation $A_n(r)$ for $n=100$, $r=1$ and $m=5, 10, 15$. We are unable to compute these values using \eqref{EQ17}, so instead we provide approximate values of these coefficients using \eqref{EQ32} and \eqref{EQ33}.

\FloatBarrier
\begin{table}[ht]
	\caption*{\textbf{Table 3.} Approximation for $A_{100}(1)$ with various values of $m$, usuing \eqref{EQ32} and \eqref{EQ33}.} 
	\label{tab:a}	
	\centering 
	\begin{tabular}{|c|c|c|} 
		\hline 
		Approximate values & Formula \eqref{EQ32} & Formula \eqref{EQ33} \\ [0.5ex] 
		\hline
		$m =5$ &  & \\
		$A_{100}(1)$ & $4.6190837447230307217\cdot10^{211}$ & $4.6190837447230037053\cdot10^{211}$\\ [0.5ex] 
		\hline
	    $m=10$ &  & \\
		$A_{100}(1)$ & $4.6190837447230307228\cdot10^{211}$ & $4.6190837447230307228\cdot10^{211}$\\[0.5ex] 
		\hline
		$m=15$ &  & \\
		$A_{40}(1)$ & $4.6190837447230307228\cdot10^{211}$ & $4.6190837447230307228\cdot10^{211}$\\ [1ex] 
		\hline
	\end{tabular}
\end{table}
\FloatBarrier


\begin{thebibliography}{99}
		\bibitem{VAdam} V.S. Adamchik, On the Barnes function. \emph{Proceedings of the 2001 International Symposium on Symbolic and Algebraic Computation}, (July 22-25, 2001, London, Canada),ACM, Academic Press, 2001, pp.15--20.
		
		\bibitem{VAdam1} Symbolic and Numeric Computations of the Barnes Function,
		Computer Physics Communications, 157(2004), 181-190.
		\bibitem{B1} E.W Barnes, The genesis of the double gamma functions, \emph{Proc. Lond. Math. Soc.} \textbf{31} (1899), pp. 358--381.
		
		\bibitem{B2} EW. Barnes, The theory of the G-function, \emph{Q. J. Pure Appl. Math.} \textbf{31} (1900), pp. 264--314.
		
		\bibitem{B3} EW. Barnes,  The theory of the double gamma function, \emph{Phil. Trans. R. Soc. Lond. A} \textbf{196}  (1901), pp. 265--388. 
		
		\bibitem{B4} EW. Barnes,  On the theory of the multiple gamma function. \emph{Trans. Camb. Philos. Soc.} \textbf{19} (1904), pp. 374--425.
		
		\bibitem{EA} E. A. Bender, An asymptotic expansion for the coefficients of some formal power series, \emph{J. Lond. Math. Soc.} \textbf{9} (1975), pp. 451--458.
		
		\bibitem{EAL}  E.A. Bender, and L.B. Richmond, An asymptotic expansion for the coefficients of some formal power series II: Lagrange Inversion. \emph{Discrete Math.} \textbf{50} (1984), pp. 135--141. 
		
		\bibitem{Berry} M. V. Berry, C. J. Howls, Infinity interpreted, \emph{Phys. World} \textbf{6} (1993), pp. 35--39.
		
		\bibitem{Boyd} W. G. C. Boyd, Approximations for the late coefficients in asymptotic expansions arising in the method of steepest descents, \emph{Methods Appl. Anal.} \textbf{2} (1995), pp. 475--489.
		
		
		\bibitem{Chen1} C.-P. Chen, Asymptotic expansions for Barnes $G$-function, \emph{J. Number Theory}, \textbf{135} (2014), pp. 36--42.
		
		
		\bibitem{JunC} J. Choi, Determinant of Laplacians on $S^3$. \emph{Math. Jpn} \textbf{40}, (1994), pp. 155--166.
		
		\bibitem{Comtet} L. Comtet, \emph{Advanced Combinotorics}, Riedel, Dordrecht, 1974.
		
		
		\bibitem{FL} C. Ferreira and J. L.  L\'opez, An asymptotic expansion of the double gamma function \emph{Journal of Approximation Theory} \textbf{111}, pp. 298--314.
		
		\bibitem{Glaisher} J. W. L. Glaisher, On a numerical continued product. \emph{Messenger of Math.} \textbf{6} (1877), pp.71--76.
		
		\bibitem{AI} A. Issaka, On Ramanujan's Inverse digamma approximation, \emph{The Ramanujan Journal} \textbf{39} (2016), pp. 291--302.
		
		
		\bibitem{Keating} J.P. Keating and N.C. Snait, Random matrix theory and $\zeta(1/2+it)$. \emph{Commun.Math.Phys} \textbf{214}, (2000), 57-89.
		
         \bibitem{Kinkelin} Kinkelin, Ueber eine mit der gammafunction verwandte transcendente und deren anwendung auf die integralrechnung. \emph{J.Reine Angew.Math.} \textbf{57} (1860), pp. 122--158.
		
		\bibitem{N1} G. Nemes, On the coefficents of the asymptotic expansion of $n!$, \emph{Integer Seqs.} \textbf{13} (2010).
		
		\bibitem{N2} G. Nemes, Approximations for the higher order coefficients in an asymptotic expansion for the Gamma function, \emph{Journal of Mathematical Analysis and Applications} \textbf{396} (2012), pp. 417--424.
		
		\bibitem{N3} G. Nemes, New asymptotic expansion for the Gamma function, \emph{Archiv der Mathematik} \textbf{95} (2010), pp. 161--169.
		
		\bibitem{N4} G. Nemes, Error bounds and exponential improvement for the asymptotic expansion of the Barnes $G$-function, \emph{Proceedings of the Royal Society A: Mathematical, Physical and Engineering Sciences.} \textbf{470} (2014), no. 2172, 14 pp.
		
		\bibitem{AMYor} A. Nikeghbali, M. Yor, The Barnes G function and its relations with sums and products of generalized Gamma, \emph{Elect. Comm. in Probab.}, \textbf{14} (2009), pp. 396--411.
		
		\bibitem{Nogu} Cassou-Nogu\'es P. 1979 Analogues p-adiques des fonctions $\Gamma$-multiples. In \emph{`Journ\'ees Arithm\'etiques de Luminy, Colloq. Internat. CNRS, Centre Univ. Luminy, Luminy, 1978'} Ast\'erisque, \textbf{61}, (1978), pp. 43--55.
		
		\bibitem{Odlyzko} A.M. Odlyzko, \emph{Asymptotic enumeration methods}. In: Graham, R.L., Gr\"tschel, M., Lov\'asz, L. (eds.)
		Handbook of Combinatorics, vol. II, pp. 1063--1229. MIT Press and North-Holland, Cambridge,
		Amsterdam (1995).
		
		
		\bibitem{OPS} B. Osgood, R. Phillips, P. Sarnak, Extremal of determinants of Laplacians. \emph{J. Funct. Anal.}
		\textbf{80}, (1988), pp. 148--211.
		
		\bibitem{QChoi} JR. Quine, J. Choi Zeta regularized products and functional determinants on spheres.
		\emph{Rocky Mountain J. Math.} \textbf{26}, (1996), pp. 719--729. 
		
		
		\bibitem{PSar} P. Sarnak, Determinants of Laplacians. \emph{Commun. Math. Phys.} \textbf{110}, (1987), pp. 113--120. 
		%
		\bibitem{PG} P.G. Todorov Taylor expansion of analytic functions related to $(1+z)^x-1$, \emph{Journal of Mathematical Analysis and Applications} \textbf{132} (1988), pp. 264--280.
		
		\bibitem{VardiI}  I. Vardi, Determinants of Laplacians and multiple gamma functions. \emph{SIAM J. Math. Anal.} \textbf{19}, (1988), pp. 493--507. 
		
		\bibitem{XW} Z. Xu and W. Wang, More asymptotic expansions for the Barnes $G$-function. \emph{Journal of Number Theory} \textbf{174} 2017, pp. 505--517.
		
		\bibitem{AminXu} A. Xu, Asymptotic expansions related to the Glaisher--Kinkelin constant and its analogues. \emph{J. Number Theory}, \textbf{163} 2016.
		
	\end{thebibliography}
\end{document}